\documentclass[11pt]{article}
\usepackage[T1]{fontenc}
\usepackage{lmodern,microtype}
\usepackage[margin=0.9in]{geometry}
\usepackage{amsmath,amssymb,amsthm,mathtools}
\usepackage[colorlinks=true,linkcolor=blue,citecolor=blue,urlcolor=blue]{hyperref}
\newcommand{\R}{\mathbb R}
\newcommand{\X}{\mathcal X}
\newcommand{\ip}[2]{\langle #1,#2\rangle}
\newcommand{\norm}[1]{\lVert #1\rVert}
\DeclareMathOperator*{\argmin}{arg\,min}
\newtheorem{theorem}{Theorem}
\newtheorem{lemma}[theorem]{Lemma}
\newtheorem{proposition}[theorem]{Proposition}
\allowdisplaybreaks[1]
\title{Projection-Free Non-Smooth Convex Programming}
\author{Kamiar Asgari \and Michael J. Neely}
\date{}
\begin{document}
\maketitle
\begin{abstract}
We study nonsmooth convex minimization over a compact convex set without
projections. For a globally $L$-Lipschitz objective, optimistic
primal--dual updates give error at most $2LD/\sqrt T$ after $T$ iterations,
where $D$ is the set's diameter. Each iteration uses at most one linear
optimization call and one subgradient query at a point that may lie outside
the feasible set. For objectives that are $\mu$-strongly convex on the
feasible set, we present a direct method with error at most
$32\rho^2L^2/[\mu(T+1)]$. It uses a local linear optimization oracle with
expansion factor $\rho$ and subgradients of norm at most $L$ queried only
at feasible points. The method requires neither restarts nor advance
knowledge of $T$. We give self-contained proofs and a short restart argument for
the convex variant of the Garber--Hazan algorithm that removes the
logarithmic factor from their error bound.
\end{abstract}

\section{Problem setup and earlier work}
We consider $\min_{x\in\X}f(x)$ for a continuous convex function $f$ and a nonempty
compact convex set $\X\subset\R^d$. Write $f^*=\min_\X f$ and
$D=\operatorname{diam}(\X)$. All norms and inner products are Euclidean.
An ordinary linear optimization oracle (LOO) is a routine that returns
$\operatorname{LOO}_\X(c)\in\argmin_{x\in\X}\ip{c}{x}$.

Asgari and Neely~\cite[Algorithm 1 and Theorem 1]{asgari2024} obtain
$2LD/\sqrt T$ for problems with only a set constraint.
Section~\ref{sec:general} gives an optimistic primal--dual interpretation
of their approach and a self-contained proof. After $T$ iterations,
the $2LD/\sqrt T$ error bound matches worst-case lower bounds of order
$LD/\sqrt T$ up to constants in sufficiently high dimensions relative to $T$,
in the respective
first-order~\cite{nemirovsky1983} and linear-oracle~\cite[Theorem~2]{lan2014}
models.

For strongly convex objectives, Section~\ref{sec:strong} gives a direct
method with error at most $32\rho^2L^2/[\mu(T+1)]$, feasible subgradient
queries, and a local linear optimization oracle (LLOO).
Section~\ref{sec:comparison} gives error at most
$1217\rho^2L^2/[\mu(T+1)]$ by restarting the Garber--Hazan algorithm.
Their existing strongly convex method has an
$O(\log(T+1)/T)$ bound~\cite[Section~6.2]{garber2015}.

\section{A primal--dual method for convex objectives}
\label{sec:general}
Assume that $f:\R^d\to\R$ is globally $L$-Lipschitz and that a subgradient
$g(y)\in\partial f(y)$ is available at every $y\in\R^d$.
Then $\|g(y)\|\le L$. If $L=0$ or $D=0$, any feasible point is optimal,
so assume $L,D>0$.

\subsection{Equivalent formulations}
To avoid projecting onto $\X$, we use two points: a feasible point
$x\in\X$ for linear optimization and an unrestricted point $y\in\R^d$
for subgradient queries. An unrestricted dual vector $q\in\R^d$ enforces $x=y$:
\begin{equation}\label{eq:general-saddle}
 f^*=\min_{\substack{x\in\X\\y\in\R^d}}
       \sup_{q\in\R^d}\bigl\{f(y)+\langle q,x-y\rangle\bigr\}.
\end{equation}
The supremum equals $f(x)$ when $x=y$ and is $+\infty$ otherwise.

Alternatively, restricting the dual vector to $\|q\|\le L$ gives
\[
 \max_{\|q\|\le L}\bigl\{f(y)+\langle q,x-y\rangle\bigr\}
 =f(y)+L\|x-y\|\ge f(x),
\]
where the inequality follows from Lipschitz continuity and is an equality
at $y=x$. Thus
\[
 f^*=\min_{\substack{x\in\X\\y\in\R^d}}
       \max_{\|q\|\le L}\bigl\{f(y)+\langle q,x-y\rangle\bigr\}.
\]
We use the unrestricted formulation \eqref{eq:general-saddle} to avoid projecting the dual
iterates onto the ball $\{q:\|q\|\le L\}$.

\subsection{Optimistic primal and dual updates}
Choose an integer $T\ge1$ and parameters $\alpha,\beta>0$.
Initialize $y_1=x_1\in\X$ and $q_1=0$.
For $t=1,\ldots,T$, query $g_t=g(y_t)\in\partial f(y_t)$ and,
if $t<T$, perform the following updates in order.
\begin{enumerate}
\item \textbf{Primal update.} Define the concave quadratic model
\[
 h_t(y)=\left\langle q_t-\frac1\beta
                  \left(\frac12y-x_t\right),y\right\rangle
\]
and set
\begin{equation}\label{eq:general-primal-update}
 y_{t+1}=\argmin_{y\in\R^d}
 \left\{\langle g_t,y\rangle-h_t(y)
                  +\frac\alpha2\|y-y_t\|^2\right\}.
\end{equation}
\item \textbf{Dual update.} With the predictions
$\tilde g_t=y_t-x_1$ for $1\le t\le T$, set
\begin{equation}\label{eq:general-dual-update}
 q_{t+1}=\argmin_{q\in\R^d}
 \left\{\left\langle q,y_t-x_t+\tilde g_{t+1}-\tilde g_t\right\rangle
                  +\frac\beta2\|q-q_t\|^2\right\}.
\end{equation}
\item \textbf{Linear optimization.} Choose
\[
 x_{t+1}=\operatorname{LOO}_\X(q_{t+1})
          \in\argmin_{x\in\X}\langle q_{t+1},x\rangle.
\]
\end{enumerate}
Return the feasible average $\bar x_T=T^{-1}\sum_{t=1}^T x_t$.
The method uses $T$ subgradient queries and $T-1$ LOO calls.
Both $y_t$ and $q_t$ are unrestricted; in particular, subgradient queries
may lie outside $\X$.

The dual minimization gives
$q_{t+1}=q_t-(y_{t+1}-x_t)/\beta$. Consequently,
\[
 \nabla h_t(y_{t+1})=q_t-\frac{y_{t+1}-x_t}{\beta}=q_{t+1}.
\]
Thus the primal optimality condition and the dual update are
\[
 \alpha(y_{t+1}-y_t)+g_t-q_{t+1}=0,\qquad
 \beta(q_{t+1}-q_t)+y_{t+1}-x_t=0.
\]
These conditions give the explicit implementation
\begin{equation}\label{eq:general-updates}
\begin{aligned}
 y_{t+1}&=\frac{\alpha\beta y_t+\beta(q_t-g_t)+x_t}{1+\alpha\beta},\\
 q_{t+1}&=q_t-\frac{y_{t+1}-x_t}{\beta},\qquad
 x_{t+1}=\operatorname{LOO}_\X(q_{t+1}).
\end{aligned}
\end{equation}
In the optimistic interpretation~\cite{rakhlin2013}, the primal and dual
descent vectors are $g_t-q_t$ and $y_t-x_t$. Using $-q_t$
and $\tilde g_t=y_t-x_1$ as predictions gives errors $g_t$ and $x_1-x_t$,
whose norms are bounded by $L$ and $D$, respectively.

\begin{theorem}\label{thm:general-rate}
For every $T\ge1$ and $\alpha,\beta>0$, the method
\eqref{eq:general-updates} satisfies
\begin{equation}\label{eq:general-bound}
 f(\bar x_T)-f^*
 \le\frac{\alpha D^2+\beta L^2}{2T}
       +\frac{L^2}{2\alpha}+\frac{D^2}{2\beta}.
\end{equation}
In particular, $\alpha=L\sqrt T/D$ and $\beta=D\sqrt T/L$ give
\[
 f(\bar x_T)-f^*\le\frac{2LD}{\sqrt T}.
\]
\end{theorem}
\begin{proof}
We first prove the prediction-error inequality used for both updates.
For $c>0$ and vectors $a_t,m_t$ with $m_1=0$, suppose
\begin{equation}\label{eq:optimistic-template}
 w_{t+1}=w_t-\frac{a_t+m_{t+1}-m_t}{c},\qquad t<T.
\end{equation}
Set $z_t=w_t+m_t/c$. Then $z_1=w_1$ and
$z_{t+1}=z_t-a_t/c$ for $t<T$. Define
$z_{T+1}=z_T-a_T/c$ only for the proof. For every comparison vector $w$,
\begin{align*}
 \langle a_t,w_t-w\rangle
 &=\langle a_t,z_t-w\rangle-\frac{\langle a_t,m_t\rangle}{c}\\
 &=\frac c2\bigl(\|z_t-w\|^2-\|z_{t+1}-w\|^2\bigr)
   +\frac{\|a_t-m_t\|^2-\|m_t\|^2}{2c}.
\end{align*}
Summing and dropping the nonpositive terms gives
\begin{equation}\label{eq:prediction-error}
 \sum_{t=1}^T\langle a_t,w_t-w\rangle
 \le\frac c2\|w_1-w\|^2
       +\frac1{2c}\sum_{t=1}^T\|a_t-m_t\|^2.
\end{equation}

Fix a minimizer $x^*$ and write $Q=\{q\in\R^d:\|q\|\le L\}$.
The primal update has the form \eqref{eq:optimistic-template} with
$w_t=y_t$, $c=\alpha$, $a_t=g_t-q_t$, and $m_t=-q_t$.
Since $m_1=0$ and $a_t-m_t=g_t$, \eqref{eq:prediction-error} yields
\begin{equation}\label{eq:general-primal-sum}
 S_y:=\sum_{t=1}^T\langle g_t-q_t,y_t-x^*\rangle
 \le\frac{\alpha D^2}{2}+\frac{TL^2}{2\alpha}.
\end{equation}
For the dual update, take $w_t=q_t$, $c=\beta$,
$a_t=y_t-x_t$, and $m_t=\tilde g_t=y_t-x_1$. Again $m_1=0$, and the
prediction error is $x_1-x_t$, so for every $q\in Q$,
\begin{equation}\label{eq:general-dual-sum}
 S_q(q):=\sum_{t=1}^T\langle y_t-x_t,q_t-q\rangle
 \le\frac{\beta L^2}{2}+\frac{TD^2}{2\beta}.
\end{equation}
The set $Q$ restricts only the comparison vector, not the dual iterates.

To combine these estimates, set
$B_T=(\alpha D^2+\beta L^2)/2+TL^2/(2\alpha)+TD^2/(2\beta)$.
Since $q_1=0$ and $x_t=\operatorname{LOO}_\X(q_t)$ for $t\ge2$,
$\langle q_t,x^*-x_t\rangle\ge0$ for every $t$.
Together with the subgradient inequality, this gives, for every $q\in Q$,
\begin{align*}
 B_T\ge S_y+S_q(q)
 &=\sum_{t=1}^T\bigl[\langle g_t,y_t-x^*\rangle
   +\langle q,x_t-y_t\rangle+\langle q_t,x^*-x_t\rangle\bigr]\\
 &\ge\sum_{t=1}^T\bigl[f(y_t)-f^*+\langle q,x_t-y_t\rangle\bigr].
\end{align*}
Define $\bar y_T=T^{-1}\sum_{t=1}^T y_t$ only for the analysis.
Averaging, maximizing over $Q$, and using convexity and Lipschitz
continuity give
\begin{equation}\label{eq:general-conversion}
 f(\bar x_T)-f^*
 \le f(\bar y_T)-f^*+L\|\bar x_T-\bar y_T\|
 \le\frac{B_T}{T}.
\end{equation}
This proves \eqref{eq:general-bound}. Substituting the stated choices
of $\alpha$ and $\beta$ gives $2LD/\sqrt T$.
\end{proof}

\section{Strong convexity: a direct method}
\label{sec:strong}
We now consider the strongly convex case, where $f$ need only be defined on $\X$.
Assume that $\mu,L>0$ are known and that each query at $x\in\X$ returns
$g(x)$ satisfying
\begin{equation}\label{eq:strong-oracle}
 \norm{g(x)}\le L,\qquad
 f(u)\ge f(x)+\ip{g(x)}{u-x}+\frac\mu2\norm{u-x}^2
 \quad(u\in\X).
\end{equation}
We also assume a local linear optimization oracle (LLOO) with known
expansion factor $\rho\ge1$. On input $z\in\X$, $r>0$, and $c\in\R^d$,
it returns $p=\operatorname{LLOO}(z,r,c)$ such that
\begin{equation}\label{eq:local-oracle}
 p\in\X,\quad \norm{p-z}\le\rho r,\quad
 \ip{c}{p}\le\ip{c}{u}\quad
 (u\in\X,\ \norm{u-z}\le r).
\end{equation}
This oracle is stronger than an ordinary LOO. Garber and Hazan construct one for polytopes~\cite[Section~5]{garber2015}.

The unique minimizer $x^*$ satisfies
\begin{equation}\label{eq:basic-strong}
 \frac\mu2\norm{x-x^*}^2\le f(x)-f^*
 \le\frac{L^2}{2\mu}.
\end{equation}
The first inequality follows by applying strong convexity between $x^*$
and $x$ and using optimality along the segment. The second follows by maximizing
the right-hand side of $f(x)-f^*\le L\norm{x-x^*}-\mu\norm{x-x^*}^2/2$
over $\norm{x-x^*}$.

\paragraph{Algorithm.}
Start from any $x_1\in\X$ and $b_0=0$. Set
$a=1/(2\rho^2)$ and, for $t=1,2,\ldots$, compute
\begin{equation}\label{eq:direct-alg}
\begin{aligned}
 A_t&=\frac{t(t+1)}2,&
 r_t&=\frac{16\rho^2L}{\mu\sqrt{t(t+1)}},& g_t&=g(x_t),\\
 b_t&=b_{t-1}+t(g_t-\mu x_t),&
 p_t&=\operatorname{LLOO}(x_t,r_t,\mu A_tx_t+b_t),\\
 x_{t+1}&=(1-a)x_t+ap_t,&
 \widehat x_T&=\frac1{A_T}\sum_{t=1}^Tt x_t.
\end{aligned}
\end{equation}
Each iteration uses one subgradient query and one LLOO call. All
iterates, query points, and outputs are feasible. The method uses no
function values and needs no prescribed number of steps or target accuracy.

\paragraph{Two model inequalities.}
We use two model inequalities in the proofs here and in Section~\ref{sec:comparison}.
Let $F_t=F_0+\sum_{i=1}^tw_i\ell_i$, $w_i>0$, and let $u_t$ minimize
$F_t$ over $\X$. The minima $m_t=F_t(u_t)$ satisfy
$m_t-m_{t-1}\ge w_t\ell_t(u_t)$. Summing this inequality gives, for any $u\in\X$,
\begin{equation}\label{eq:model-telescope}
 \sum_{t=1}^Tw_t[\ell_t(x_t)-\ell_t(u)]
 \le F_0(u)-\min_\X F_0+
       \sum_{t=1}^Tw_t[\ell_t(x_t)-\ell_t(u_t)].
\end{equation}
If $F_t$ is quadratic with Hessian $\kappa_t I$, $\kappa_t>0$, write
$d_t=F_t(x_t)-m_t$ and $e_t=F_t(x_{t+1})-m_t$ for the model errors
before and after the update.
For an LLOO step with cost vector $\nabla F_t(x_t)$ and
stepsize $a\in[0,1]$, if $\norm{x_t-u_t}\le r_t$, a quadratic expansion gives
\begin{equation}\label{eq:local-descent}
 e_t\le(1-a)d_t+\frac{\kappa_t}{2}a^2\rho^2r_t^2.
\end{equation}
Indeed, $\ip{\nabla F_t(x_t)}{p_t-x_t}
\le\ip{\nabla F_t(x_t)}{u_t-x_t}\le-d_t$.
The model minimizers $u_t$ are used only in the proof.

\begin{theorem}\label{thm:direct}
The method \eqref{eq:direct-alg} satisfies, for every $T\ge1$,
\begin{equation}\label{eq:direct-rate}
 f(\widehat x_T)-f^*\le\frac{32\rho^2L^2}{\mu(T+1)}.
\end{equation}
\end{theorem}
\begin{proof}
Take $F_0=0$, $w_t=t$, and
$\ell_t(x)=\ip{g_t}{x}+\mu\norm{x-x_t}^2/2$.
Then $\nabla^2F_t=\mu A_t I$ and
$\nabla F_t(x_t)=\mu A_tx_t+b_t$.
Set
\[
 M=\frac{64\rho^4L^2}{\mu},\qquad
 \chi=1-\frac1{4\rho^2},\qquad e_0=0.
\]
We first prove $d_t\le M$ and $e_t\le\chi M$ by induction; these bounds
also justify the choice of radius $r_t$.
Suppose $e_{t-1}\le\chi M$. Comparing successive minima and using strong
convexity of $F_t$ gives
\[
 d_t\le e_{t-1}+t[\ell_t(x_t)-\ell_t(u_t)]
 \le\chi M+tL\norm{x_t-u_t}
 \le\chi M+\frac{2L}{\sqrt\mu}\sqrt{d_t}.
\]
This also holds at $t=1$ since $F_0=0$.
Since $(1-\chi)M=2L\sqrt{M/\mu}$, solving the quadratic inequality in
$\sqrt{d_t}$ gives $d_t\le M$. Thus
$\norm{x_t-u_t}\le\sqrt{2M/(\mu A_t)}=r_t$ and
\eqref{eq:local-descent} gives
$e_t\le(1-a+a^2\rho^2)M=\chi M$. This completes the induction.

It remains to translate the model bounds into a function-value bound.
Strong convexity, averaging, and \eqref{eq:model-telescope} give
\[
\begin{aligned}
 A_T[f(\widehat x_T)-f^*]
 &\le\sum_{t=1}^Tt[\ell_t(x_t)-\ell_t(x^*)]\\
 &\le\sum_{t=1}^Tt[\ell_t(x_t)-\ell_t(u_t)]
 \le\sum_{t=1}^TtLr_t
 \le\frac{16\rho^2L^2T}{\mu}.
\end{aligned}
\]
Dividing by $A_T=T(T+1)/2$ proves \eqref{eq:direct-rate}.
\end{proof}

The weights $t$ keep $tLr_t$ bounded while $A_T$ grows quadratically,
giving the $32\rho^2L^2/[\mu(T+1)]$ bound. Such weights also
appear in projected dual averaging~\cite[Corollary~1]{tao2021}; here one LLOO step
per iteration controls the model error.

\section{Restarting the Garber--Hazan algorithm}
\label{sec:comparison}
Under the assumptions of Section~\ref{sec:strong}, restarting the
Garber--Hazan algorithm gives error at most $1217\rho^2L^2/[\mu(T+1)]$.
We first prove a bound for one stage, then choose a restart schedule
that halves the error bound at each stage.

\paragraph{One stage.}
For a starting point $z$, stepsize $\eta>0$, and integer $N\ge1$,
define the procedure $\mathcal C(z,\eta,N)$ as follows. Set
$x_1=z$, $K=1+18\rho^2$, and $a=1/(3\rho^2)$. For $t=1,\ldots,N$,
query $g_t=g(x_t)$ and compute
\begin{equation}\label{eq:gh-convex}
\begin{aligned}
 F_t(x)&=\norm{x-z}^2+\eta\sum_{i=1}^t\ip{g_i}{x},&
 r&=K\eta L,\\
 p_t&=\operatorname{LLOO}(x_t,r,\nabla F_t(x_t)),&
 x_{t+1}&=(1-a)x_t+ap_t,\\
 \mathcal C(z,\eta,N)&=\frac1N\sum_{t=1}^Nx_t.
\end{aligned}
\end{equation}
These updates follow the $\sigma=0$ branch of the Garber--Hazan
algorithm~\cite[Algorithm~5 and Lemma~11]{garber2015}, with $\eta$ left free.

\begin{lemma}\label{lem:gh-convex}
For every $z,u\in\X$, the procedure satisfies
\begin{equation}\label{eq:gh-convex-rate}
 f(\mathcal C(z,\eta,N))-f(u)
 \le\frac{\norm{z-u}^2}{\eta N}+K\eta L^2.
\end{equation}
\end{lemma}
\begin{proof}
Set $F_0(x)=\norm{x-z}^2$, let $u_t$ minimize $F_t$ over $\X$, and use
the model errors $d_t,e_t$ defined in Section~\ref{sec:strong}, with
$u_0=z$ and $e_0=0$.
Write $h=(18\rho^2\eta L)^2$. Since the Hessian of $F_t$ is $2I$,
\[
 \norm{u_t-u_{t-1}}^2
 \le F_t(u_{t-1})-F_t(u_t)
 \le\eta\ip{g_t}{u_{t-1}-u_t}
 \le\eta L\norm{u_t-u_{t-1}}.
\]
Thus $\norm{u_t-u_{t-1}}\le\eta L$.
Suppose $e_{t-1}\le h$. Then
$\norm{x_t-u_t}\le\sqrt h+\eta L=r$ and
$d_t\le h+\eta L\sqrt h+\eta^2L^2$.
Equation~\eqref{eq:local-descent} gives
\[
 e_t\le(1-a)(h+\eta L\sqrt h+\eta^2L^2)+a^2\rho^2r^2
 =h-\left(54\rho^2+1+\frac{2}{9\rho^2}\right)\eta^2L^2\le h.
\]
The induction starts at $e_0=0$. Apply \eqref{eq:model-telescope} with
$w_t=\eta$ and $\ell_t(x)=\ip{g_t}{x}$, and use convexity:
\[
 N[f(\mathcal C(z,\eta,N))-f(u)]
 \le\frac{\norm{z-u}^2}{\eta}
       +\sum_{t=1}^N\ip{g_t}{x_t-u_t}
 \le\frac{\norm{z-u}^2}{\eta}+NLr.
\]
Divide by $N$.
\end{proof}

\paragraph{Restart schedule.}
Each stage halves the stepsize and doubles the number of steps.
Set $z_0=z$,
$m=\lceil64K\rceil$, and, for $j=0,1,\ldots$, run
\begin{equation}\label{eq:restart}
 \Delta_j=\frac{L^2}{2\mu\,2^j},\qquad
 \eta_j=\frac{\Delta_j}{4KL^2},\qquad N_j=m2^j,
 \qquad z_{j+1}=\mathcal C(z_j,\eta_j,N_j).
\end{equation}
Each stage resets the model on the same set $\X$ with the same LLOO.
After an integer $T\ge0$ steps, return $z_J$, where
$T_J=m(2^J-1)\le T<T_{J+1}$. Thus $z_J$ is the last completed stage output,
or $z_0$ if $T<m$.

\begin{proposition}\label{prop:restart}
The restart schedule \eqref{eq:restart} satisfies
\begin{equation}\label{eq:restart-rate}
 f(z_J)-f^*\le\frac{mL^2}{\mu(T+m)}
 \le\frac{1217\rho^2L^2}{\mu(T+1)}.
\end{equation}
\end{proposition}
\begin{proof}
Equation~\eqref{eq:basic-strong} gives $f(z_0)-f^*\le\Delta_0$.
If $f(z_j)-f^*\le\Delta_j$, then
$\norm{z_j-x^*}^2\le2\Delta_j/\mu$ and
$N_j\ge32KL^2/(\mu\Delta_j)$. Lemma~\ref{lem:gh-convex} gives
\[
 f(z_{j+1})-f^*
 \le\frac{2\Delta_j}{\mu\eta_jN_j}+K\eta_jL^2
 \le\frac{\Delta_j}{4}+\frac{\Delta_j}{4}=\Delta_{j+1}.
\]
Thus, at each stage boundary,
$f(z_J)-f^*\le mL^2/[2\mu(T_J+m)]$.
To cover steps in an unfinished stage, use $T<T_{J+1}$, which gives
$T+m<2m2^J$, so
$\Delta_J\le mL^2/[\mu(T+m)]$.
Finally $m\le64(1+18\rho^2)+1\le1217\rho^2$.
\end{proof}

The direct and restarted methods each use at most $T$ subgradient
queries in $\X$ and $T$ LLOO calls, with no advance choice of $T$.
The direct method has the smaller proved constant and returns an updated
average at every step; the restarted method returns the last completed
stage's average.

\section*{AI assistance}
During the preparation of this work, the authors used GPT-6 to assist with developing and presenting parts of the analysis for strongly convex objectives.

{\small
\bibliographystyle{plain}
\bibliography{references}
}
\end{document}